\documentclass[12pt,leqno]{article}
\usepackage{amsmath,amssymb,amsthm,amscd,euscript,mathrsfs}
\input xy
\xyoption{all}

\newtheorem{lemma}{Lemma}

\newtheorem{corollary}{Corollary}
\newtheorem{definition}{Definition}
\newtheorem{proposition}{Proposition}
\newtheorem{theorem}{Theorem}
\newtheorem{remark}{Remark}

\usepackage{verbatim}
\usepackage{fancyhdr}

 \def\mP{\mathbb P}

 \def\kO{\mathcal O}

\begin{document}

\title{Curves on a nonsingular Del Pezzo Surface in $\mP^4_k$}
\author{Elena Drozd}
\date{}

\maketitle

\begin{abstract}

Classification of curves in a projective space occupies minds of many mathematicians.
First step in doing so is classification of curves on a given surface. This brings us to 
consideration of the nonsingular Del Pezzo Surface in $\mP^4_k.$ We describe conditions
for a divisor to be very ample and  condition for a divisor class to contain an irreducible curve.

\end{abstract}


~\newline
Study of the surfaces and curves that lie on them is a subject treated by many authors.
Various types of surfaces were considered. 
Mori \cite{mori} studied smooth curves on smooth 
quadric surface in $\mP^3$.  A.Knutsen  \cite{knutsen} studied   
smooth curves on $K3$ surfaces in $\mP^3$ .  
Nonsingular cubic surfaces in $\mP^3_k$ and the famous 27 lines on them were studied
by R.Hartshorne in \cite{hartshorne}.  
Following the work of R. Hartshorne \cite {hartshorne}  we study the existence of 
nondegenerate integral curves on a nonsingular Del Pezzo surfaces in $\mP^4_k$ by 
representing Del Pezzo Surface in $\mP^4_k$ as a $\mP^2_k$ with 5 points blown up.
We describe conditions for a divisor to be very ample and  condition for a divisor 
class to contain an irreducible curve.

Let $k$ be an algebraically closed field of characteristic $0$.
\begin{definition}\label{def21}
A nonsingular \emph{Del Pezzo surface} is defined to be a surface $X$ of degree $4$ in $\mP^4_k$
 such that $w_X \cong {\mathcal O}_X(-1)$.
\end{definition}

Now we recall a theorem which will give us the construction of the nonsingular Del Pezzo 
surfaces.

\begin{theorem}[\cite{hartshorne} V.4.6]\label{theorem:hartshorne1}
Let $S$ be the linear system of plane cubic curves with assigned (ordinary) base points
$P_1,\dots,P_r$, and assume that no $3$ of the $P_i$ are collinear, and no 6 of them lie on a conic.
If $r\leq 6$, then the correspondingly linear system $S^\prime$ on the surface
$X$ obtained from $\mP^2_k$ by blowing up $P_1,\dots,P_r$, is very ample.
\end{theorem}

\begin{corollary}[\cite{hartshorne} V.4.7]\label{corollary:hartshorne2}
With the same hypotheses as in \ref{theorem:hartshorne1}, for each $r=0,1,\dots,6$ we obtain 
an embedding of $X$ in $\mP^{9-r}_k$ as a surface of degree $9-r$, whose canonical sheaf $w_X$ is 
isomorphic to ${\mathcal O}_X(-1)$. In particular, for $r=5$, we obtain a Del Pezzo surface of
degree $4$ in $\mP^r_k$.
\end{corollary}

In fact, every nonsingular Del Pezzo surface in $\mP^4_k$ can be obtained in this way 
(\cite{manin} 24.4). Also, any Del Pezzo surface is a complete intersection of two nonsingular 
quadric hypersurfaces. 

To prove this, let $X$ be a Del Pezzo surface. We need to find $h^0\bigl(\kO_X(2)\bigr)$. 
By the Rieman-Roch theorem for surfaces we obtain:
$$h^0(D)-h^1(D)+h^2(D)=\frac12 D.(D-K)+1+p_a,$$
where $p_a$ is arithmetic genus of $X$ and $p_a=0$ since $X$ is a rational surface. 
Let $D=2H$. Then
$$\frac12 D.(D-K)+1+p_a=\frac12(2H). (3H)+1=3H^2+1=13$$
since $K=-H$ and $H^2=4=\deg X$. By \cite[7.7]{hartshorne} \quad $H^2\bigl(X,\kO_X(2)\bigr)$ is 
dual to 
$H^0\bigl(\kO_X(-3)\bigr)$, and therefore it is zero.

Now, we show that $h^1\bigl(X,\kO_X(2)\bigr)=0$. \quad $\kO_X(1)$ corresponds to a divisor $H$ 
which is a nonsingular elliptic quartic curve $C$. From the exact sequence
$$0\to \kO_X \to \kO_X(H) \to \kO_H(H) \to 0$$ 
taking cohomology we obtain
$$\dots \to H^1(\kO_X)\to H^1\bigl(\kO_X(H)\bigl)\to H^1\bigl(\kO_H(H)\bigr)\to \dots \quad . 
$$
$h^1\bigl(\mP^2_k,\kO_{\mP^2_k})=0$, therefore by blowing up theorem \cite{hartshorne} V.3.4 \quad 
$h^1(X,\kO_X)=0$. Also, $h^1\bigl(\kO_H(H)\bigr)=0$ since $\kO_H(n), ~ n\geq 1$ is a nonspecial divisor for
a degree 4 and genus 1 nonsingular curve. Therefore 
\begin{equation}\label{equation:h1zero}
h^1\bigl(\kO_X(H)\bigr)=0
\end{equation}
Now in the exact sequence 
$$0\to \kO_X(H)\to \kO_X(2H)\to \kO_H(2H)\to 0$$ 
taking cohomology we obtain that $h^1\bigl(\kO_X(2H)\bigr)=0$ since $h^1\bigl(\kO_H(2H)\bigr)=0$ as 
$h^1\bigl(\kO_H(H)\bigr)=0$, and $h^1\bigl(\kO_X(H)\bigr)=0$ by (\ref{equation:h1zero}).
Thus we arrive a $h^0\bigl(X,\kO_X(2H)\bigr)=13$, while $h^1\bigl(\kO_{\mP^4_k}(2)\bigr)=15$.
Thus $X$ lies on two quadric hypersurfaces $Q_1$ and $Q_2$. Since $\deg X=\deg Q_1\cdot \deg Q_2$ 
we must have $X=Q_1\cap Q_2$. Therefore $X$ is a complete intersection of $Q_1$ and $Q_2$.

Conversely, 
a nonsingular surface that is a complete intersection of two quadric hypersurfaces is a Del Pezzo 
surface.
Suppose $X$ is a nonsingular surface such that $X=Q_1\cap Q_2$ is a complete intersection.
Then $\deg X=\deg Q_1\cdot \deg Q_2=4$ and by \cite{hartshorne} ex.II.8.4 (e)
$$w_X=\kO_X(\deg Q_1 +\deg Q_2 -5) = \kO_X(-1).$$
Therefore $X$ is a Del Pezzo surface by \ref{def21}.

To study curves on the nonsingular Del Pezzo surface, we establish some notation.

Let $P_1,\dots,P_5$ be five points of $\mP^2_k$, no three collinear (we can ommit condition ``no 
6 of them lie on a cone" since we have only 5 points).

Let $S$ be the linear system of plane cubic curves through $P_1,\dots,P_5$, and let $X$ be the 
nonsingular Del Pezzo surface in $\mP^4_k$ obtained by the \ref{theorem:hartshorne1}, and 
\ref{corollary:hartshorne2}. Thus $X$ is isomorphic to $\mP^2_k$ with 5 points
$P_1,\dots,P_5$ blown up. Let $\pi : X \longrightarrow \mP^2_k$ be the projection. 
Let $E_1,\dots,E_5\subseteq X$ be the exceptional curves, and let $e_1,\dots,e_5 \in 
\text{\textrm{Pic }} X$ be their linear equivalence classes. Let $l\in \textrm{Pic } X$ be the 
class of $\pi^*$ (of a line in $\mP^2_k$).

\begin{proposition}\label{proposition:prop1}
Let $X$ be the Del Pezzo surface in $\mP^4_k$ constructed as above. Then:
\begin{enumerate}
\item\label{item:item1} $\text{\emph{Pic }} X \cong {\mathbb Z}^6$, generated by $l,e_1,\dots,e_5$;
\item\label{item:item2} the intersection pairing on $X$ is given by 
$l^2=1,e_i^2=-1,l\cdot e_i=0,e_i\cdot e_j = 0 \text{ for } i\neq j$;
\item\label{item:item3} the hyperplane section $h$ is $3l-\sum_{i=1}^5 e_i$;
\item\label{item:item4} the canonical class is $K=-h=-3l+\sum_{i=1}^5 e_i$;
\item\label{item:item5} if $D$ is any effective divisor on $X$, $D\sim al-\sum_{i=1}^5 b_i e_i$ 
for some $a\geq 0$, then the degree of $D$, as a curve in $\mP^4_k$ is $d=3a-\sum_{i=1}^5 b_i$;
\item\label{item:item6} the self-intersection of $D$ is $D^2=a^2-\sum_{i=1}^5 b_i^2$;
\item\label{item:item7} the arithmetic genus of $D$ is 
$$p_a(D)=\frac{1}{2}(a-1)(a-2)-\frac{1}{2}\sum_{i=1}^5 b_i (b_i-1)=\frac{1}{2}(D^2-d)+1 \text{.}$$
\end{enumerate}
\begin{proof}
Cases (\ref{item:item1}) and (\ref{item:item2}) follow from (\cite{hartshorne}, V.3.2.). 
Case (\ref{item:item3}) comes from the definition of the 
embedding in $\mP^4_k$. Case (\ref{item:item4}) comes from (\cite{hartshorne}, V.3.3.). 
For case (\ref{item:item5}), note that degree 
of $D$ is just $D\cdot h$. Case (\ref{item:item6}) follows from case (\ref{item:item2}). 
Case (\ref{item:item7}) follows from the adjunction formula 
$2p_a(D)-2=D\cdot (D+K)$ and the fact that $D\cdot K=-D\cdot H=-d$, by (\ref{item:item4}).
\end{proof}
\end{proposition}

\begin{remark}\label{remark:rem1}
If $C$ is any irreducible curve in $X$, other then $E_1,\dots,E_5$, then $\pi (C)$ is an 
irreducible plane curve $C_0$, and $C$ is the strict transform of $C_0$. Let $C_0$ have degree $a$. 
Suppose that $C_0$ has a point of multiplicity $b_i$ at each $P_i$. Then
$\pi^*C_0=C+\sum_{i=1}^5b_iE_i$ by (\cite{hartshorne}, V.3.6). Since $C_0\sim a\cdot \text{(line) }$ 
we conclude that $C\sim al-\sum_{i=1}^5b_ie_i$ . Thus for any $a,b_1,\dots,b_5\geq 0$, we can 
interpret an irreducible curve $C$ on $X$ in the class $al-\sum_{i=1}^5b_ie_i$ as the strict 
transform of a plane curve of degree $a$ with a $b_i$-fold point at each $P_i$. So the study of 
curves on $X$ is reduced to the study of certain plane curves.
\end{remark}

\begin{theorem}[Lines on Del Pezzo surface]\label{theorem:lines}
The Del Pezzo surface in $\mP^4_k$ contains exactly 16 lines. Each one has self-intersection $-1$, 
and they are the only irreducible curves with negative self-intersection on $X$. They are 
\begin{enumerate}
\item the exeptional curves $E_i$ ($i=1,\dots,5$) (five of these),
\item the strict transforms $F_{ij}$ of the line in $\mP^2_k$ containing $P_i$ and $P_j$, 
$1 \leq i < j \leq 5$ (ten of these), and
\item the strict transform $G$ of the conic in $\mP^2_k$ containing all five points $P_1,\dots,P_5$.
\end{enumerate}
\begin{proof}
If $L$ is any line in $X$, then $\text{deg} L=1$, $p_a(L)=0$, so by \ref{proposition:prop1} part (7) 
we have $0=\frac{1}{2} (L^2-1)+1$, therefore $L^2=-1$. Conversely, if $C$ is an irreducible curve on 
$X$ with $C^2<0$, then since $p_a(C)\geq 0$, we must have $\frac{1}{2} (C^2-d)+1\geq 0$ or 
$C^2\geq d-2$. Thus $C^2=-1$, $p_a(C)=0$, $\text{deg} C=1$, so $C$ is a line.

Next, $E_i\sim e_i$, $F_{ij}\sim l-e_i-e_j$, $G\sim 2l-\sum_{i=1}^5 e_i$, and each of these has 
degree $1$, $p_a=0$  (from \ref{proposition:prop1}) i.e. is a line.

It remains to show that if $C$ is any irreducible curve on $X$ with $\text{deg} C=1$, and $C^2=-1$, 
then $C$ is one of those 16 lines listed.

Assume $C$ is not one of the $E_i$. Then we can write $C\sim al-\sum_{i=1}^5b_ie_i$ and by Remark 
\ref{remark:rem1} we must have $a>0, b\geq 0$. Furthermore, 
\begin{equation*}
\text{deg} C=3a-\sum_{i=1}^5b_i=1\text{,} \qquad
C^2=a^2-\sum_{i=1}^5b_i^2=-1 \text{.}
\end{equation*}
We will show that the only integers $a,b_i,\dots ,b_5$ satisfying all these conditions are those 
corresponding to $F_{ij}$ and $G$ above.

By Schwartz's inequality: if $x_1,\dots,x_n,y_1,\dots,y_n$ are two sequences of real numbers, then 
$$\mid \sum_{i=1}^nx_iy_i\mid ^2\leq \mid \sum_{i=1}^nx_i^2 \mid \cdot \mid  \sum_{i=1}^ny_i^2 \mid 
\text{.}$$
Taking $x_i=1,i=1,\dots,5, y_i=b_i, i=1,\dots,5$ we get 
$$\bigl(\sum_{i=1}^5b_i\bigr)^2\leq 5\bigl(\sum_{i=1}^5b_i^2\bigr)\text{.}$$
But from the conditions above $$\sum_{i=1}^5b_i=3a-1, \sum_{i=1}^5b_i^2=a^2+1 
\text{, therefore we obtain}$$
$$(3a-1)^2\leq 5(a^2+1)\text{\quad or \quad} 4a^2-6a-4<0.$$
Solving  this inequality we get $a\leq 2$. Therefore, $a=1$ or $a=2$. 
Now one finds all possible values of the $b_i$ by trial:
if $a=1$, then $b_i=b_j=1$ for some $i\neq j$ and the rest $0$. 
This gives $F_{ij}$. If $a=2$ then all $b_i=1$. This gives $G$.
\end{proof}
\end{theorem}

\begin{proposition}\label{proposition:prop2}
Let $X$ be a Del Pezzo surface as above, and let $E^\prime_1,\dots,E^\prime_5$ be any subset of 
five mutually skew lines chosen from among 16 lines on $X$. Then there is another morphism 
$\pi^\prime:X\longrightarrow \mP^2_k$, making $X$ isomorphic to that $\mP^2_k$ with five points 
$P_1^\prime,\dots,P_5^\prime$ blown up (no 3 collinear), such that $E^\prime_1,\dots,E^\prime_5$ 
are the exceptional curves for $\pi^\prime$.
\begin{proof}
We will show first that it is possible to find $\pi^\prime$ such that $E_1^\prime$ is the inverse 
image of $P^\prime_1$.

\underline{Case 1.}  If $E^\prime_1$ is one of the $E_i$, we take $\pi^\prime=\pi$, but relabel 
$P_i$ so that $P_i$ becomes $P_1^\prime$.

\underline{Case 2.} If $E^\prime_1$ is one of the $F_{ij}$, say $E^\prime=F_{12}$, then we apply 
the quadratic transformation with centers $P_1, P_2, P_3$
(see \cite{hartshorne} V.4.2.3) as follows. Let $X_0$ be $\mP^2_k$ with $P_1, P_2, P_3$ blown up, let	 
$\pi_0: X_0\longrightarrow \mP^2_k$ be the projection, and let  $\psi: X_0\longrightarrow \mP^2_k$ be
the other map to$\mP^2_k$ (see picture below), so that $X_0$ via $\psi$ is 
$\mP^2_k$ with $Q_1, Q_2, Q_3$ blown up.

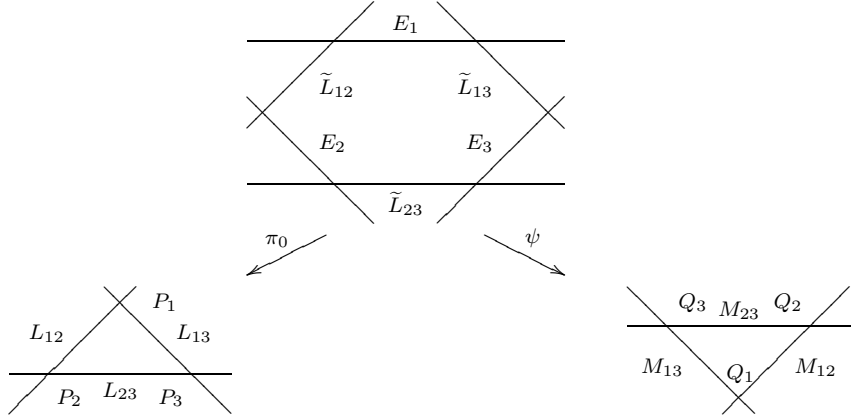
\begin{equation*}
\xymatrix@=12pt{
& & & & & & & & & {} \ar@{-}[lllddd]^{ \widetilde{L}_{12}} & \ar@{-}[rrrddd]_{\widetilde{L}_{13}} & & & & & & & & &\\
& & & & & & {} \ar@{-}[rrrrrrr]^{E_1} & & & & & & &  & & & & & & \\
& & & & & & {} \ar@{-}[rrrddd]^{E_2} & & & & & & & \ar@{-}[lllddd]_{E_3} & & & & & &\\
& & & & & &   & & & & & &  & & & & & & & \\
& & & & & & {} \ar@{-}[rrrrrrr]_{\widetilde{L}_{23}} & & & & & & & & & & & & &  \\
& & & & \ar@{}[ddd]^{P_1} & &  &  & \ar[lld]_{\pi_0}  & & &  \ar[rrd]^{\psi} & & & & & & & \\
& & & \ar@{-}[dddrrr]^{L_{13}} & \ar@{-}[dddlll]_{L_{12}} & & & & &   &        & & &  & \ar@{-}[rrrddd]_{M_{13}} & & & & &  \ar@{-}[lllddd]^{M_{12}} \\	
& & & & & & & & &                                                                                        &     & & &  &  \ar@{-}[rrrrr]^{M_{23}} & & & & & \\	
&\ar@{-}[rrrrr]_{L_{23}} & \ar@{}[d]^{P_2} & & & \ar@{}[d]_{P_3} & & & & & & & &  & & &  \ar@{}[uuu]^{Q_3}  & & \ar@{}[uuu]^{Q_2} & \\	
& & & & & & & & &                                                                                              & & & &  & & & &  \ar@{}[uu]^{Q_1} & & \\	
& & & & & & & & &                                                                                              & & & &  & & & & & & \\	
}
\end{equation*}

\underline{note}: $\psi^{-1}(Q_i)=\widetilde{L}_{jk}$

Since $\pi: X_0\longrightarrow \mP^2_k$ expresses $X$ as $P^2$ with $P_1,\dots,P_5$ blown up, $\pi$ factors
through $\pi_0$, say $\pi=\pi_0\circ \theta$ where $\theta : X\longrightarrow X_0$. Now we define $\pi^\prime$ as $\psi \circ \theta$.

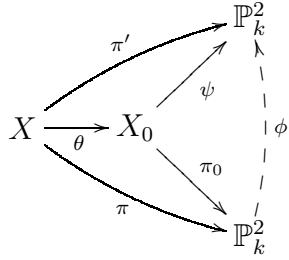
\begin{equation*}
\xymatrix{
 & & \mP^2_k \\
X \ar@/^/[urr]^{\pi^\prime} \ar[r]_\theta \ar@/_/[drr]_\pi & X_0 \ar[ur]_\psi \ar[dr]^{\pi_0} & \\
 & & \mP^2_k \ar@/_/@{-->}[uu]_\phi
}
\end{equation*}
Then $\theta(F_{12})=\widetilde{L}_{12}$, so $\pi^\prime(F_{12})=Q_3$. Furthermore, $\pi^\prime$ expresses $X$ as $\mP^2_k$ with
$Q_1$, $Q_2$, $Q_3$, $P_4^\prime$, $P_5^\prime$ blown up, where $P_4^\prime,P_5^\prime$  are the images of $P_4$ and $P_5$ under $\psi\circ\pi_0^{-1}$.
Now, taking $P_1^\prime=Q_3$, and $P_2^\prime, P_3^\prime$ to be $Q_1,Q_2$, we have $E_1^\prime={\pi^\prime}^{-1}(P_1^\prime)$.

We still have to verify that no 3 of $Q_1, Q_2, Q_3, P_4^\prime, P_5^\prime$ lie on a line. $Q_1, Q_2, Q_3$ are non-collinear by construction.

If $Q_1, Q_2, P_4^\prime$ were collinear, then $\psi^{-1}(P_4^\prime)\in E_3$, so $P_4$ would be infinitely near $P_3$, but we assumed that 
$P_1,\dots,P_5$ are ordinary points.

If $Q_1, P_4^\prime, P_5^\prime$ were collinear, then let $L^\prime$ be the line containing them. Consider its strict transform $S$ by $\phi^{-1}$.
It is a line $L$ containing $P_1, P_4, P_5$ (contradiction to no three of $P_1,\dots,P_5$ collinear). Indeed, $\phi^{-1}$ is the rational
map determined by the linear system of conics through $Q_1, Q_2, Q_3$. Such a conic has one free intersection with $L^\prime$,
so the strict transform of $L^\prime$ is a line $L$. Furthermore, $L^\prime$ meets $M_{23}$, so $L$ passes through $P_1$. This
completes case 2.

\underline{Case 3.}  If $E_1^\prime$ is $G$. Apply the quadratic transformation with centers $P_1$, $P_2$, $P_3$. 

Since $\pi(G)$ is the conic through $P_1$, $P_2,\dots,P_5$ we see that $\pi^\prime(G)$ is the line through $P_4^\prime$, 
$P_5^\prime$. Thus $E_1^\prime$ is the curve $F_{45}^\prime$ for $\pi^\prime$, which reduces to the case 2.

Now that we have moved $E_1^\prime$ to the position of $E_1$, we may assume that $E_1^\prime=E_1$, and
consider $E_2^\prime$. Since $E_2^\prime$ doesn't meet $E_1$, the possible values of $E_2^\prime$ are $E_2$, $E_3$, $E_4$, $E_5$,
$F_{ij}$ with $1<i,j$. We apply the same method as in cases 1,2 above, and find that we can move $E_2^\prime$ to the 
role of $E_2$ without changing $P_1$. That is to say, we allow to relable only $P_2,\dots,P_5$, or use quadratic transformation based at
three points among $P_2,\dots,P_5$.

As in case 1 - case 2 we can move $E_2^\prime$ to $E_2$ without changing $P$.

For $E_3^\prime$ possible values are $E_3, E_4, E_5, F_{34}$, $F_{35}$, $F_{45}$.

So, by either relabling points  $P_3$, $P_4$, $P_5$ or using quadratic transformation based at points $P_3$, $P_4$, $P_5$ we
can move $E_3^\prime$ to $E_3^\prime$. 

Now for $E_4$ possible values are $E_4, E_5$, $F_{45}$, but $F_45$ meets $E_4$ and $E_5$, therefore lines
$E_4^\prime$, $E_5^\prime$ must be $E_4$ and $E_5$ in some order. So, for the last step we have only to permute 4 and 5 if
necessary.

\end{proof}
\end{proposition}

\begin{remark}\label{remark28}
The proposition says that any five mutually skew lines among the 16 lines play the role of $E_1,\dots,E_5$.
\end{remark}

\begin{theorem}
The following conditions are equivalent for a divisor $D$ on the Del Pezzo surface $X$:
\begin{enumerate}
\item $D$ is very ample;
\item $D$ is ample;
\item $D^2>0$, and for every line $L\subseteq X$, $D.L>0$;
\item for every line $L\subseteq X$, $D.L>0$.
\end{enumerate}
\begin{proof}

To prove this theorem we need the following:

\begin{lemma}\label{lemma:l310}
Let $D\sim al-\sum_{i=1}^5b_ie_i$ be a divisor class on the Del Pezzo surface $X$, and suppose that $b_1\geq b_2\geq\dots\geq b_5>0$
and $a\geq b_1+b_2+b_5$. Then $D$ is very ample.
\begin{proof}
We use the general fact that a very ample divisor plus a divisor moving in a linear system without base points is very ample.

Let us consider the divisor classes

\begin{align*}
D_0 &= l \\
D_1 &= l-e_1 \\
D_2 &= 2l-e_1-e_2 \\
D_3 &= 2l- e_1 - e_2 - e_3 \\
D_4 &= 2l - e_1 - e_2 - e_3 - e_4 \\
D_5 &= 3l - e_1 - e_2 - e_3 - e_4 - e_5 
\end{align*}

Then $|D_0|$, $|D_1|$ correspond to the linear system of lines in $P^2$ with 0 or 1 assigned base points, which have no
unassigned base points. $|D_2|$, $|D_3|$, $|D_4|$ have no basis points by (\cite{hartshorne}, V.4.1.).
Therefore any linear combination of these $D=\sum_{i-1}^5c_iD_i$ with $c_i\geq 0$, $c_5>0$ will be very ample divisior.

Clearly, $D_0,\dots,D_5$ form a free basis for $\text{Pic} X\cong {\mathbb Z}^6$. Wrting $D\sim al-\sum_{i=1}^5b_ie_i$, we
have 
\begin{multline*}
c_0l+c_1(l-e_1)+ c_2 (2l-e_1-e_2 ) +
c_3 ( 2l- e_1 - e_2 - e_3 ) + \\
+ c_4 ( 2l - e_1 - e_2 - e_3 - e_4 )+
c_5 ( 3l - e_1 - e_2 - e_3 - e_4 - e_5) 
= al-\sum_{i=1}^5b_ie_i \text{,}
\end{multline*}
\begin{align*}
\text{therefore \qquad{}} a &=  c_0 + c_1+2c_2+2c_3+2c_4+3c_5 \\
b_5 &= c_5\\
b_4 &= c_5+c_4 \\
b_3 &= c_5+c_4 + c_3 \\
b_2 &= c_5+c_4 + c_3 +c_2 \\
b_1 &= c_5+c_4 + c_3 +c_2 +c_1 \text{\qquad with $c_i\geq 0, c_5>0$, therefore }
\end{align*}
\begin{multline*}
a = c_0+(c_1+ c_2  + c_3 +c_4+ c_5) + (c_2  + c_3 +c_4+ c_5) + c_5 \text{,} \\
\text{therefore \quad} 
a=c_0+b_1+b_2+b_5 \text{,\quad therefore \quad } a\geq b_1+b_2+b_5 \text{.} 
\end{multline*}
Thus, conditions $c_i\geq 0, c_5>0$ are equivalent to the conditions $b_1\geq b_2 \geq \dots \geq b_5 >0$, 
$a\geq b_1+b_2+b_5$, therefore all divisors satisfying these conditions are very ample.
\end{proof}
\end{lemma}
Proof of the Theorem:

Of course (1) $\Rightarrow$ (2) $\Rightarrow$ (3) $\Rightarrow$  (4) using easy direction of Nakai's criterion
(\cite{hartshorne}, V.1.10). Now for (4) $\Rightarrow$ (1): Suppose $D$ is a divisor satisfying $D.L>0$ for every line $L\subseteq X$.
Choose five mutually skew lines $E_1^\prime,\dots,E_5^\prime$ as follows:

Choose $E_5^\prime$ so that $D.E_5^\prime$ is equal to the minimum value of $D.L$ for any line $L$;

Choose $E_4^\prime$ so that $D.E_4^\prime$ is equal to the minimum value of $D.L$ among those lines $L$ which 
do not meet $E_5^\prime$.

Choose $E_3^\prime$ so that $D_1E_3^\prime$ is equal to the minimum value of $D.L$ among those lines $L$ which 
do not meet $E_4^\prime$, $E_5^\prime$.

Now there are only three lines left which do not meet $E_3^\prime$, $E_4^\prime$, $E_5^\prime$. One of them meeting the other two.
Choose $E_2^\prime$, $E_1^\prime$ so that $D.E_1^\prime\geq D.E_2^\prime$.

Now by \ref{remark28} we may assume that $E_i^\prime=E_i$, $i=1,\dots,5$. Writing $D\sim al-\sum_{i=1}^5b_ie_i$ we have
$D.E_i=b_i$, so by construction we have $b_1\geq b_2\geq \dots \geq b_5 > 0$. On the other hand, $F_{12}$ was available
as a candidate at the time we chose $E_3$ so we must have $D.F_{12}\geq D.E_3$.

This translates as $a-b_1-b_2\geq b_3$, i.e. $a\geq b_1+b_2+b_3$, but $b_3\geq b_5$, therefore $a\geq b_1+b_2+b_5$,
therefore these conditions imply that $D$ is very ample.
\end{proof}
\end{theorem}
\begin{corollary}\label{corollary:delPezzoIrr}
A divisor class $D$ on the Del Pezzo surface contains an irreducible curve $\Leftrightarrow$ it contains an irreducible
nonsingular curve $\Leftrightarrow$ it is either (a) one of the 16 lines, or (b) a conic with $D^2=0$, or (c) $D.L\geq 0$ for 
every line $L$ and $D^2>0$.
\begin{proof}
Let $X_i$ be a surface obtained from $\mP^2_k$ by blowing up points $P_1,\dots,P_i, \quad i\leq 4$.
Then one can prove analog of lemma \ref{lemma:l310} for such surfaces. If in the divisor class 
$D=(a; b_1,\dots b_i)$ on Del Pezzo surface 
\begin{enumerate}
\item $b_5=0, b_4\neq 0$, consider $X_4$,
\item $b_5=b_4=0, b_3\neq 0$, consider $X_3$,
\item $b_5=b_4=b_3=0, b_2\neq 0$, consider $X_2$,
\item $b_5=b_4=b_3=b_2=0, b_3\neq 0$, consider $X_1$,
\item $b_5=b_4=b_3=b_2=b_1=0, b_3\neq 0$, consider $\mP^2_k$.
\end{enumerate}
In each of the cases above divisor class $(a; b_1,\dots,b_i)$ has a
nonsingular irreducible curve, say, $C_i$, which, if we rechoose points $P_{i+1},\dots,P_5$
so that $C_i$ does not pass through $P_{i+1},\dots,P_5$, will correspond to a nonsingular
irreducible curve in the divisor class $D$. This concludes the proof of the corollary.
\end{proof}
\end{corollary}

\newpage

\end{document}